\documentclass[11pt]{amsart}
\usepackage{amsmath,amsthm, amscd, amssymb, amsfonts, mathrsfs,pb-diagram,pb-xy}
\usepackage[all]{xy}
\usepackage{enumitem}
\usepackage{mathrsfs}
\usepackage{srcltx}
\usepackage{hyperref}
\usepackage{multicol}
\usepackage[dvips, dvipsnames, usenames]{color}

\usepackage{ textcomp }

\allowdisplaybreaks

\newcommand{\superf}{{\mathbf F}(4)}
\newcommand{\superg}{{\mathbf G}(3)}

\newcommand{\superqa}[3]{{\bf A}_{#1}(#2|#3)}
\newcommand{\superqb}[3]{{\bf B}_{#1}(#2|#3)}

\newcommand{\superqd}[3]{{\bf D}_{#1}(#2|#3)}

\newcommand{\adc}{\operatorname{ad}_c}
\newcommand{\co}{\operatorname{co}}

\newcommand{\gr}{\operatorname{gr}}

\newcommand{\id}{\operatorname{id}} 
\newcommand{\ord}{\operatorname{ord}}

\newcommand{\GK}{\operatorname{GKdim}}

\newcommand{\ku}{ \Bbbk}
\newcommand{\I}{\mathbb I}
\newcommand{\G}{\mathbb G}
\newcommand{\N}{\mathbb N}
\newcommand{\Z}{\mathbb Z}
\newcommand{\J}{\mathbb J}

\newcommand{\toba}{\mathscr{B}}
\newcommand{\wtoba}{\widetilde{\mathscr{B}}}
\newcommand{\htoba}{\widehat{\toba}}
\newcommand{\bq}{\mathfrak{q}}

\newcommand{\ydhc}{{}^{H_{0}}_{H_{0}}\mathcal{YD}}

\newcommand{\Pc}{\mathcal{P}}

\newcommand{\Zc}{\mathcal{Z}}

\def\br{\mathfrak{br}}

\def\bgl{\mathfrak{wk}}

\def\g{\mathfrak{g}}

\newcommand{\pre}{\mathfrak{Pre}(V)}
\newcommand{\prefd}{\mathfrak{Pre}_{\operatorname{fGK}}(V)}

\newcommand{\qti}{\widetilde{q}}

\newcommand{\pf}{\begin{proof}}
\newcommand{\epf}{\end{proof}}

\numberwithin{equation}{section}
\theoremstyle{plain}

\newtheorem{theorem}{Theorem}[section]
\newtheorem{lemma}[theorem]{Lemma}
\newtheorem{definition-theorem}[theorem]{Definition-Theorem}
\newtheorem{prop}[theorem]{Proposition}

\newtheorem{corollary}[theorem]{Corollary}

\theoremstyle{definition}

\newtheorem{remark}[theorem]{Remark}

\newcommand{\ot}{\otimes}

\begin{document}
\noindent
\title[Finite GK-dimensional pre-Nichols algebras]
{Pre-Nichols algebras of one-parameter families of finite Gelfand-Kirillov dimension}

\author[E. Campagnolo]{Emiliano Campagnolo}
\email{emiliano.campagnolo@mi.unc.edu.ar}

\keywords{Hopf algebras, Nichols algebras, Gelfand-Kirillov dimension.
\\
MSC2020: 16T05, 16T20, 17B37, 17B62.}

\thanks{The work was supported by CONICET}

\begin{abstract}
Between the braided vector spaces of diagonal type there exist some families whose associated Nichols algebras are infinite dimensional but with finite Gelfand-Kirillov dimension, called one-parameter families. 
We show that every pre-Nichols algebra of finite Gelfand-Kirillov dimension of this kind with connected diagram is necessarily the corresponding Nichols algebra, up to few exceptions. For each one of these exceptions we present a proper pre-Nichols algebra by generators and relations, and prove that it is the corresponding eminent pre-Nichols algebra.

This paper essentially ends the determination of the eminent pre-Nichols algebras of diagonal type, 
the first step towards the classification of pre-Nichols algebras of diagonal type with finite GKdim. 
\end{abstract}

\maketitle

\section{Introduction}

Let $\Bbbk$ be an algebraically closed field of characteristic zero. If $H$ is a pointed Hopf algebra (that is, its coradical $H_0$ coincides with the subalgebra generated by the group $G(H)$ of group-likes elements) over $\Bbbk$ then the coradical filtration $\{H_i\}_{i\in \N_0}$ of $H$ is a Hopf algebra filtration with $H_0= \ku G(H)$, thus the associated graded object $\gr H$ is a graded Hopf algebra and admits morphisms of Hopf algebras $\pi:\gr H\rightarrow \ku G(H)$ and $\iota:\ku G(H)\rightarrow \gr H$ with $\pi\circ\iota=\id$. By \cite[\S 11.7]{Rad-libro}, $\gr H\simeq \ku G(H)\# R$, where $R=\bigoplus_{n\ge 0} R(n)$ is a graded Hopf algebra in the Yetter-Drinfeld modules category $\ydhc$ such that $R(0)=\Bbbk 1$. Let $V:= R(1)\in \ydhc$: the subalgebra of $R$ generated by $V$ is the Nichols algebra $\toba(V)$ of $V$, which and is a graded Hopf algebra in $\ydhc$ univocally determined by $V$.

Given a group $G$, a program (called the \emph{Lifting Method}) for the classification of finite dimensional pointed Hopf algebras $H$ such that $G(H)\simeq G$ is presented in \cite{AS}. There is a similar program when the Gelfand-Kirillov dimension is finite, which proposes to answering the following questions:
\begin{enumerate}[leftmargin=*]
\item\label{lifting-m-1} Classify all Yetter-Drinfeld modules $V$ over $G$ such that $\GK \toba(V)<\infty$.
\item\label{lifting-m-2} For each $V$ obtained above, compute all post-Nichols algebras of $V$ with $\GK$ finite.
\item Find all Hopf algebras $H$ such that $\gr H\simeq \mathcal{E}\#\Bbbk G$.
\end{enumerate}
From the main result in \cite{angiono-garcia-2022} it is proved that every Nichols algebra of diagonal type with finite $\GK$ has a finite root system. On the other hand in \cite{H-classif} all braided vector spaces $V$ with connected Dynkin diagram and finite root systems were classified. Thus, combining these two results we answer \eqref{lifting-m-1} when $V$ is of diagonal type; that is, we obtain the classification of all Nichols algebras of diagonal type with connected Dynkin diagram and finite $\GK$.

The problem proposed in \eqref{lifting-m-2} for diagonal type is equivalent to computing all pre-Nichols algebras with $\GK$ finite. 
The main results of \cite{ASa,ACSa,ACSa2} give an answer for \eqref{lifting-m-2} when the Nichols algebra is finite-dimensional: they prove the existence of a pre-Nichols algebra $\htoba$ of $V$ such that every other finite $GK$-dimensional pre-Nichols algebra of $V$ is a quotient of $\htoba$; that is, $\htoba$ is the corresponding eminent pre-Nichols algebra.

This paper gives a complete answer for \eqref{lifting-m-2} when $V$ is a braided vector space of diagonal type such that $\dim \toba(V)=\infty$ but $\GK \toba(V)<\infty$; these braided vector spaces are obtained by evaluating the parameter of the \emph{one-parameter families} \cite{AAY} in a generic value. The strategy follows the same ideas as in the series of papers \cite{ASa,ACSa,ACSa2} but here we describe the whole poset of pre-Nichols algebras with finite $\GK$, which has just one or two elements as we describe now:

\begin{theorem}\label{thm:introduction-infinite-parameters}
Let $V$ be a finite-dimensional braided vector space of diagonal type with connected Dynkin diagram such that $\dim \toba(V)=\infty$ but $\GK \toba(V)<\infty$.
\begin{enumerate}[leftmargin=*,label=\rm{(\roman*)}]
\item\label{item:introduction-infinite-parameters-i} If $V$ is not of type $\superqa{3}{q}{\{2\}}$, $\superqa{3}{q}{\{1,2,3\}}$ or $D(2,1;\alpha)$ such that exactly two of the parameters $q,r,s$ are not roots of unity, then $\toba(V)$ is the only pre-Nichols algebra of $V$ with $\GK<\infty$.
\item\label{item:introduction-infinite-parameters-ii} If $V$ is of type $\superqa{3}{q}{\{2\}}$, $\superqa{3}{q}{\{1,2,3\}}$ or $D(2,1;\alpha)$ with diagram \eqref{diagram:first-exceptional-case-D21a}, \eqref{diagram:second-exceptional-case-D21a}, \eqref{diagram:third-exceptional-case-D21a}, then the set of finite $\GK$ pre-Nichols algebras has two elements: $\toba(V)$ and $\htoba(V)$, where $\htoba(V)$ is defined by \eqref{eq:superA-first-case-eminent}, respectively \eqref{eq:superA-second-case-eminent}, \eqref{eq:D(2,1,a)-first-case}, \eqref{eq:D(2,1,a)-second-case}, \eqref{eq:D(2,1,a)-third-case}.
\end{enumerate}
\end{theorem}

Indeed, if $V$ is as in \ref{item:introduction-infinite-parameters-i}, then $\toba(V)$ is the unique pre-Nichols algebra with finite $\GK$ by Theorem \ref{thm:nichols-is-eminent}.

If $V$ is as in \ref{item:introduction-infinite-parameters-ii}, then $\htoba(V)$ is eminent by Propositions \ref{prop:superA-first-case-eminent}, \ref{prop:superA-second-case-eminent}, \ref{prop:D(2,1,a)-first-case}, \ref{prop:D(2,1,a)-second-case}, \ref{prop:D(2,1,a)-third-case}. Notice that 
$\Pc(\htoba(V))=V\oplus \Bbbk\{z\}$, for $z$ as in Remarks \ref{rem:caseA1}, \ref{rem:caseA2}, \ref{rem:caseD1}, \ref{rem:caseD2}, \ref{rem:caseD3}, hence the unique Hopf quotient of $\htoba(V)$ that preserves $V$ is $\toba(V)=\htoba(V)/\langle z \rangle$.

A difference of these results with those in \cite{ASa,ACSa,ACSa2}, is that here one finds all pre-Nichols algebras with finite $\GK$ while for the case  $\dim\toba(V)<\infty$ it remains to determine the whole poset of pre-Nichols algebras with finite $\GK$. 

The structure of the paper is the following. Section \ref{sec:preliminaries} contains the basic theory of diagonal type Nichols algebras required for the rest of the work. In Section \ref{sec:lemmas} we prove a series of lemmas on pre-Nichols algebras of diagonal type analogous to those of \cite{ACSa,ACSa2} that serve for the proof of the main results. These results are contained in Section \ref{sec:Eminent-infinite}, where we present the eminent pre-Nichols algebra for each one of the cases. In the same section we consider the non connected case: we are able to determine the whole poset when each connected component is as in Theorem \ref{thm:introduction-infinite-parameters}.   
%

\section{Preliminaries}\label{sec:preliminaries}

In this section we recall general aspects and particular results of the theory of finite Gelfand-Kirillov dimensional pre-Nichols algebras of diagonal type. For the basic theory and notation on Nichols algebras we refer to \cite{A-leyva}.

\subsection*{Notations}
The set of $N$-th roots of unity of $\ku^{\times}$ is denoted by $\G_N$. Also $\G_{\infty}:=\displaystyle \cup_{N\geq 0}\G_N$. 

Let $\theta\in\N$. We denote by $\I_{\theta}=\{1,\dots ,\theta\}$, and simply write $\I$ if $\theta$ is fixed. 
The element $(a_1,a_2,\cdots,a_{\theta})\in\Z^{\I}$ is denoted by $1^{a_1}2^{a_{2}}\cdots \theta^{a_{\theta}}$.
For each $i\in\I$, let $\alpha_i=i^1$, the $i$-th element of the canonical basis of $\Z^{\I}$.

For notation related to Hopf algebras we refer to \cite{Rad-libro}, and for the definition and basic facts of Gelfand-Kirillov dimension see \cite{KL}.

\subsection{Eminent pre-Nichols algebras}

Let $V$ be a finite-dimensional braided vector space. A pre-Nichols algebra of $V$ is a Hopf quotient $\toba$ of $T(V)$ by an $\N_0$-homogeneous ideal generated by elements of degree $\ge 2$; hence, there exist canonical epimorphisms $T(V) \twoheadrightarrow \toba \twoheadrightarrow \toba(V)$, induced by the identity of $V$. 

Let $\pre$ be the poset whose elements are the pre-Nichols algebras of $V$ and whose partial order is given by $\toba_1\leq\toba_2$ if there exists an epimorphism $\toba_1\twoheadrightarrow\toba_2$ induced by the identity on $V$. Therefore $T(V)$ is minimum in $\pre$ and $\toba(V)$ is maximum in $\pre$. 

We denote by $\prefd$ the subposet of $\pre$ that contains all pre-Nichols algebras of $V$ with finite $\GK$. Notice that $T(V)\notin \prefd$ if $\dim V>1$, and $\prefd$ is not empty if and only if $\GK\toba(V)<\infty$. If there exists a pre-Nichols algebra $\htoba$ that is minimal on $\prefd$, we will call it the \emph{eminent} pre-Nichols algebra of $V$. 

The main results of this work determine the eminent pre-Nichols algebra for some cases when $V$ is a braided vector space of diagonal type.

\subsection{Nichols algebras of diagonal type with finite Gelfand-Kirillov dimension}
We will use here the same notation for Nichols algebras of diagonal type as in \cite[\S 2]{ACSa2}, see also the definitions and results stated there. 
For example, given a braided vector space $(V,c)$ of diagonal type with braided matrix $\bq=(q_{ij})_{i,j\in\I}$, we denote by $\toba_\bq$ the associated Nichols algebra. 
Next we recall results that will be useful for subsequent sections.

\medspace

Let $(x_i)_{i\in\I}$ be a basis such that $c(x_i\ot x_j)=q_{ij} \, x_j\ot x_i$, $i,j\in\I$.
Every Nichols algebra $\toba_\bq$ is $\N_0^{\I}$-graded, where $\deg x_i=\alpha_i$.
By \cite{Kh} $\toba_\bq$ has a (truncated) PBW basis whose generators are $\N_0^{\I}$-homogeneous: the set $\varDelta^\bq_+$ of degrees of the homogeous generators does not depend on the PBW basis and is called the set of positive roots of $\toba_{\bq}$.
There is a close relation between the $\GK$ of $\toba_\bq$ and its positive roots; namely:

\begin{theorem}\label{thm:conjecture-true}\cite{angiono-garcia-2022}
$\GK\toba_\bq<\infty$ if and only if the set $\varDelta^\bq_+$ is finite.
\end{theorem}

The classification of Nichols algebras with finite root system is given in \cite{H-classif}. Thus, if $\GK\toba_\bq<\infty$, then the Dynkin diagram of $\bq$ appears in \cite[Tables 1--4]{H-classif}.
These Nichols algebras are grouped in \cite{AA17} into families according with the nature of the associated root system: Cartan, super, (super) modular and unidentified types.
Also, a presentation by generators and relations of each example is exposed in \cite{AA17} (following \cite{An-crelle}). 

In this paper we focus on those Nichols algebras from \cite{H-classif} such that $\dim \toba(V)=\infty$:
\begin{enumerate}[leftmargin=*,label=\rm{(\alph*)}]
\item Cartan types $A_{\theta}$, $B_{\theta}$, $C_{\theta}$, $D_{\theta}$, $E_{6,7,8}$, $F_4$, $G_2$ with parameter $q\not\in\G_{\infty}$,
\item super types $\superqa{\theta}{q}{\J}$, $\superqb{\theta}{q}{\J}$, $\superqd{\theta}{q}{\J}$, $\superf$ and $\superg$ with parameter $q\not\in\G_{\infty}$, and $D(2,1;\alpha)$ with parameters $q,r,s$ such that $qrs=1$ and at least one of them $\notin \G_{\infty}$,
\item modular types $\bgl(4)$, $\br(2)$ with parameter $q\not\in\G_{\infty}$.
\end{enumerate}

\medspace

The proofs in Section \ref{sec:lemmas} refer repeatedly to the classification of Nichols algebras in \cite{H-classif}. To avoid looking over the whole list every time, we present here some criteria for a matrix $\bq$ satisfying $\GK \toba_{\bq}=\infty$.

\begin{corollary}\label{cor:conditions-of-diagrams}
Let $\bq$ be such that $\GK\toba_\bq<\infty$. Then:
\begin{enumerate}[leftmargin=*]
\item\label{item:cycles}\cite[Lemma 20]{H-classif} The diagram of $\bq$ does not contain  $N$-cycles for all $N\geq 4$.
\item\label{item:rank3}\cite[Lemma 7(ii)]{H-classif} Let $i,j,k\in\I$ belong to a $3$-cycle, then $(q_{ii}+1)(q_{jj}+1)(q_{kk}+1)=0$ and $\qti_{ij}\qti_{jk}\qti_{ki}=1$. Morever if $q_{ii}=-1$ and $q_{jj},q_{kk}\ne -1$ then $q_{jj}\qti_{ij}=q_{kk}\qti_{ik}=1$. 
\item\cite[Proposition 4.16]{AAH-memoirs}\label{item:1-conected} The diagram does not contain a subdiagram of the form  $\xymatrix@C=30pt{\overset{1}{\circ} \ar@{-}[r]^{p}\ &\overset{q}{\circ}}$, $p\ne 1$.
\end{enumerate}
\end{corollary}



Let $\bq$ be such that $\dim \toba_{\bq}<\infty$. In this case there exist proper pre-Nichols algebras with finite $\GK$. 
One possibility is the distinguished pre-Nichols algebra $\wtoba_{\bq}$, introduced in \cite{An-dist}.
The determination of the eminent pre-Nichols algebra of $\bq$ is found in the main results of \cite{ASa,ACSa,ACSa2} and is related with $\wtoba_{\bq}$. Namely:

\begin{theorem}\label{thm:finite-case} If $\bq$ is not of type
\begin{enumerate}[leftmargin=*]
\item\label{item-exc-1} Cartan $A_{\theta}$ or $D_{\theta}$ with $q=-1$,
\item\label{item-exc-2} Cartan $A_2$ with $q\in \G_3'$,
\item\label{item-exc-3} $\superqa{3}{q}{\{2\}}$ or $\superqa{3}{q}{\{1,2,3\}}$, with $q\in \G_{\infty}$,
\item\label{item-exc-4} $\g(2,3)$ with any of the following Dynkin diagram
\begin{align*}
&d_1: \, \xymatrix@C=30pt{\overset{-1}{\circ} \ar@{-}[r]^{\xi^2} &\overset{\xi}{\circ}\ar@{-}[r]^{\xi}\ &\overset{-1}{\circ}}, &
&d_2: \, \xymatrix@C=30pt{\overset{-1}{\circ} \ar@{-}[r]^{\xi} &\overset{-1}{\circ}\ar@{-}[r]^{\xi}\ &\overset{-1}{\circ}},
\end{align*}
\end{enumerate}
then the distinguished pre-Nichols algebra $\wtoba_{\bq}$ is eminent.

If $\bq$ is of type (\ref{item-exc-2}), (\ref{item-exc-3}) or (\ref{item-exc-4}), then the eminent pre-Nichols algebra $\htoba_{\bq}$ is a central extension of $\wtoba_{\bq}$ explicitly described in \cite[Proposition 4.11]{ASa}, \cite[Proposition 5.5, 5.9]{ACSa} and \cite[Proposition 4.2,4.3]{ACSa2}, respectively.
\end{theorem}

\begin{remark}
Let $\bq$ be such that $\dim \toba_{\bq}<\infty$. If $x_{\beta}=0$ is a relation of the fixed presentation of $\wtoba_{\bq}$, then $x_{\beta}\ne 0$ in $\htoba_{\bq}$ if only if $q_{i\beta}=1$ for all $i\in\I$. The proof is by inspection for all cases such that $\htoba_{\bq}\ne \wtoba_{\bq}$.
\end{remark}

\section{Relations of finite GK-dimensional pre-Nichols algebra of $\toba$}\label{sec:lemmas}

Let $\bq$ be a braided matrix of diagonal type such that $\dim\toba_{\bq}=\infty$ and $\GK\toba_{\bq}<\infty$. Fix $\toba$ a pre-Nichols algebra of $\bq$ such that $\GK \toba <\infty$. In this section we prove whether a defining relation of $\toba_{\bq}$ in \cite{AA17} also vanishes in $\toba$. Each relation will be developed in a separate lemma whose proof follows the same strategy as in \cite[\S 3]{ACSa} and  \cite[\S 3.2]{ACSa2}. 

We start by analysing quantum Serre relations. We consider two cases: $m_{ij}=0$ and $m_{ij}\ne 0$.

\begin{lemma}\label{lem:mij=0}
Let $i,j\in \I_{\theta}$ be such that $m_{ij}=0$. If $\bq$ is not of type $\superqa{3}{q}{\{1,2,3\}}$ then $x_{ij}= 0$ in $\toba$.
\end{lemma}
\pf
Suppose that $x_{\beta}:=x_{ij}\neq 0$. Note that $x_{\beta}\in\Pc(\toba)$ since $x_{\beta}$ is a defining relation of minimal degree of $\toba_{\bq}$. 
Also
\begin{align*}
\qti_{i\beta}=q_{ii}^2&&\qti_{j\beta}=q_{jj}^2&& q_{\beta\beta}=q_{ii}q_{jj}
\end{align*}

By inspection, each matrix $\bq$ satisfies one of following conditions:
\begin{enumerate}[leftmargin=*,label=\rm{(\Roman*)}]
\item\label{item:mij=0:qii-or-qjj-not-root} $q_{ii}\notin\G_{\infty}$ or $q_{jj}\notin\G_{\infty}$. 
\item\label{item:mij=0:qii=qjj=-1} $q_{ii}=q_{jj}=-1$ and there exists $k\in \I_{\theta}$, $k\ne \{i,j\}$ such that $\qti_{ik}\qti_{jk}\not\in\G_{\infty}$.
\end{enumerate}
Assume that \ref{item:mij=0:qii-or-qjj-not-root} holds. Without loss of generality we assume that $q_{jj}\notin\G_{\infty}$. Since the Dynkin diagram is connected, there exists $r \geq 1$ and vertices $i_1,i_2,...,i_r$ such that $\qti_{ii_1} \not = 1 \neq \qti_{i_rj}$ and $\qti_{i_si_{s+1}}\neq 1$ for all $1\leq s\leq r-1$. Then $q_{ii}=-1$: otherwise, the subdiagram with the vertices $\beta,i,i_1,i_2,...,i_r,j$ in the Dynkin diagram of $W:=V+\Bbbk x_{\beta}$ is a $(r+3)$-cycle, which is not possible by Corollary \ref{cor:conditions-of-diagrams} (\ref{item:cycles}). Thus the Dynkin diagram of $\Bbbk x_j\oplus\Bbbk x_{\beta}$ is
\begin{align*}
\xymatrix@C=30pt{\overset{q_{jj}}{\underset{j}{\circ}}\ar@{-}[rr]^{q_{jj}^2} && \overset{-q_{jj}}{\underset{\beta}{\circ}}
}
\end{align*}
but this diagram does not belong to \cite[Table 1]{H-classif}, so $\GK \toba(W)=\infty$ by Theorem \ref{thm:conjecture-true}. On the other hand, $\GK \toba(W) \le \GK \toba$ by \cite[Lemma 2.8]{ASa} and we get a contradiction. Hence $x_{\beta}=0$ in $\toba$.

Next assume that \ref{item:mij=0:qii=qjj=-1} holds. As $\qti_{k\beta}=\qti_{ik}\qti_{jk}\ne 1$, the Dynkin diagram of $\Bbbk x_{\beta}\oplus\Bbbk x_k$ is
\begin{align*}
\xymatrix@C=30pt{\overset{q_{kk}}{\underset{k}{\circ}}\ar@{-}[rr]^{\qti_{ik}\qti_{jk}} && \overset{1}{\underset{\beta}{\circ}}
}
\end{align*}
so $x_{\beta}=0$ by Corollary \ref{cor:conditions-of-diagrams} (\ref{item:1-conected}).
\epf

\begin{lemma}\label{lem:mij>0}
Let $i,j\in\I_\theta$ be such that $m_{ij}>0$ and $q_{ii}^{m_{ij}+1}\neq 1$. Then $(\adc x_i)^{m_{ij}+1}x_j=0$ in $\toba$.
\end{lemma}
\pf
Suppose that $x_{\beta}:=(\adc x_i)^{m+1}x_j \neq 0$. As in the previous case, $x_{\beta}\in\Pc(\toba)$.
Set $q=q_{ii}$ and $m=m_{ij}$. The Dynkin diagram of $\ku x_i\oplus\ku x_j\oplus\ku x_\beta$ is:
\begin{align}\label{diagram-qs}
\xymatrix@C=30pt{\overset{{q_{jj}}}{\underset{j}{\circ}} \ar@{-}[rrrr]^{q^{-m}} \ar@{-}[drr]_{q^{-m(m+1)}q_{jj}^2}&&  && \overset{q}{\underset{i}{\circ}} \ar@{-}[dll]^{{q^{m+2}}}\\
&&\overset{{q^{m+1}q_{jj}}}{\underset{\beta}{\circ}}&&
}
\end{align}
By inspection, $q\not\in\G_{\infty}$. If $q^{-m(m+1)}q_{jj}^2\ne 1$, then \eqref{diagram-qs} is a triangle with at most one vertex labelled with $-1$ not satisfying the conditions in Corollary \ref{cor:conditions-of-diagrams} \eqref{item:rank3}, a contradiction. Thus $q^{-m(m+1)}q_{jj}^2=1$ and the Dynkin diagram of $\ku x_i\oplus\ku x_j\oplus\ku x_\beta$ is
\begin{align*}
& \xymatrix@C=30pt{\overset{q_{\beta\beta}}{\underset{\beta}{\circ}} \ar@{-}[rr]^{q^{m+2}}&&\overset{q}{\underset{i}{\circ}}   \ar@{-}[rr]^{q^{-m}}& &\overset{q_{jj}}{\underset{j}{\circ}}, &q\not\in\G_{\infty},
} 
\end{align*}
which does not belong to \cite[Table 2]{H-classif}. Therefore $x_{\beta}=0$ in $\toba$.
\epf

Next we deal with some powers of PBW generators annihilating in $\toba_{\bq}$.

\begin{lemma}\label{lem:xi-no-Cartan} 
Let $i\in\I_\theta$ be a non-Cartan vertex, $N=\ord q_{ii}$. Then $x_i^N=0$ in $\toba$. 
\end{lemma}
\pf
Analogous to \cite[Lemma 3.1]{ACSa}.
\epf

\begin{lemma}\label{lem:xij^2}
Let $i,j\in \I_{\theta}$ be such that $q_{ii}=\qti_{ij}=q_{jj}=-1$ and there exists $k\in\I_{\theta}-\{i,j\}$ such that $\qti_{ik}^2\qti_{jk}^2\neq 1$. Then $x_{ij}^2=0$.
\end{lemma}
\pf
A pair $i,j$ as above appears only for $\bgl(4)$. The proof is analogous to that of \cite[Lemma 3.17]{ACSa}.
\epf

Now we study a defining relation between generators attached to vertices of a triangle.

\begin{lemma}\label{lem:triangle-x_{ijk}}
Let $i,j,k\in \I_{\theta}$ be such that $\qti_{ij},\qti_{ik},\qti_{jk}\neq 1$. Then 
$$ x_{ijk} - q_{ij}(1-\qti_{jk})x_{j}x_{ik}+\frac{1-\qti_{jk}}{q_{kj}(1-\qti_{ik})}[x_{ik},x_{j}]_c = 0\text{ in }\toba.$$
\end{lemma}
\pf
Analogous to \cite[Lemma 3.22]{ACSa}.
\epf

The subsequent Lemmas deal with other relations appearing for some of the matrices $\bq$. We specify, in the proof of each lemma, for which families these are defining relations.

\begin{lemma}\label{lem:[[[xijk,xj],[xijkl,xj]],xjk]} 
Let $i,j,k,\ell \in \I_{\theta}$ and $q\in\Bbbk-(\G_2\cup\G_3)$ be such that 
\begin{align*}
q_{\ell \ell}&=\qti_{k\ell}^{\,-1}=q_{kk}=\qti_{jk}^{\,-1}=q^2, &
q_{jj}&=-1, & q_{ii}&=\qti_{ij}^{\,-1}=q^{-3}, & \qti_{ik}&=\qti_{i\ell}=\qti_{j\ell}=1.
\end{align*}
Then 
$[[[x_{ijk},x_j]_c,[x_{ijk\ell},x_j]_c]_c,x_{jk}]_c=0$.
\end{lemma}
\pf
A set of vertices as above only appears when $\bq$ is of type $\superf$. The proof is analogous to that of \cite[Lemma 3.23]{ACSa}.
\epf

\begin{lemma}\label{lem:diagonal-xijkxj}
Let $i,j,k\in\I_\theta$ be such that $q_{jj}=-1$, $\qti_{ik}=\qti_{ij}\qti_{jk}=1$ and $\qti_{ij} \neq \pm 1$. If $\bq$ is not of type $\superqa{3}{q}{\{2\}}$ then $[x_{ijk},x_j]_c=0$ in $\toba$.
\end{lemma}
\pf
By inspection, $\bq$ satifies one of the following conditions:
\begin{enumerate}[leftmargin=*]
\item\label{item:xijkj=0-some--1} $q_{ii}=-1$ or $q_{jj}=-1$;
\item \label{item:case-b-xijkj=0} $q_{ii}q_{kk}=1$ and there exists  $\ell\in\I_\theta-\{i, j, k\}$ such that $\qti_{i\ell}\neq 1=\qti_{j\ell}=\qti_{k\ell}\,$;
\item \label{item:case-c-xijkj=0} $q_{ii}q_{kk}=1$ and there exists $\ell\in\I_\theta-\{i, j, k\}$ such that $\qti_{j\ell}^2\neq 1=\qti_{i\ell}=\qti_{k\ell}\,$;
\item \label{item:case-d-xijkj=0} $q_{ii}q_{kk}=1$ and there exists $\ell\in\I_\theta-\{i, j, k\}$ such that $\qti_{k\ell}\neq 1=\qti_{j\ell}=\qti_{i\ell}$.
\end{enumerate}
The corresponding proofs are analogous to \cite[Lemma 3.13, 3.14]{ASa}
\epf

\begin{lemma}
Let $i,j,k\in \I_{\theta}$ be such that $q_{ii}=q_{jj}=-1$, $\qti_{ij}^{\,2}=\qti_{jk}^{\,-1}\ne 1$ and $\qti_{ik}=1$. Then $[[x_{ij},x_{ijk}]_c,x_j]_c=0$ in $\toba$.
\end{lemma}
\pf
This relation appears when $\bq$ is of type $\superqd{n}{q}{\J}$, $\superf$ or $\superg$. By inspection, either $q_{kk}^2\ne1$ or else $\qti_{ij}^{\, 3}\ne 1$. The proof is analogous to that of \cite[Lemma 3.19]{ACSa}.
\epf

\begin{lemma}
Let $i,j,k,\ell\in \I_{\theta}$ be such that $q_{kk}= -1$, $q_{jj}\qti_{ij}=q_{jj}\qti_{jk}=1$, $\qti_{ik}=\qti_{il}=\qti_{j\ell}=1$ and $\qti_{jk}^2=\qti_{k\ell}^{-1}=q_{\ell\ell}$. Then $[[[x_{ijk\ell},x_{k}]_c,x_{j}]_c,x_{k}]_c=0$ in $\toba$.
\end{lemma}
\pf
A set of vertices as above only appears when $\bq$ is of type $\superqd{n}{q}{\J}$. The proof is analogous to that of \cite[Lemma 3.21]{ACSa}.
\epf

\begin{lemma}
Let $i,j,k \in \I_{\theta}$ be such that $q_{ii}=q_{jj}=-1$, $\qti_{ij}^{\ 3}= \qti_{jk}^{\ -1}$ and $\qti_{ik}=1$. 
Then $[[x_{ij},[x_{ij},x_{ijk}]_c]_c,x_j]_c=0$ in $\toba$.
\end{lemma}
\pf
A set of vertices as above only appears when $\bq$ is of type $\superg$. The proof is analogous to that of \cite[Lemma 3.26]{ACSa}.
\epf

\begin{lemma}
Let $i,j,k \in \I_{\theta}$ be such that $q_{ii}\notin\G_{\infty}$, $\qti_{ik}=1$, $q_{jj}=-1$, $\qti_{ij}=q_{ii}^{-2}$, $\qti_{jk}=q_{kk}^{-1}=-q_{ii}^{3}$.
Then 
\begin{align*}
&[x_i,[x_{ijk},x_j]_c]_c -\frac{q_{ij}q_{kj}}{1-q_{ii}^{-1}}[x_{ij},x_{ijk}]_c -(q_{ii}+q_{ii}^{2}) q_{ij}q_{ik}x_{ijk}x_{ij}=0 && \text{ in }\toba.
\end{align*}
\end{lemma}
\pf
Again, a set of vertices as above only appears when $\bq$ is of type $\superg$. The proof is analogous to that of \cite[Lemma 3.27]{ACSa}.
\epf



\begin{lemma}\label{lem:[[xijkl,xj],xk]} Let $i,j,k,\ell \in \I_{\theta}$ be such that one of the following hold:
\begin{enumerate}[leftmargin=*,label=\rm{(\roman*)}]
\item \label{item:[[xijkl,xj],xk]-a}
$q_{kk}=-1$, $q_{ii}=\qti_{ij}^{\,-1}=q_{jj}^2$, $\qti_{k\ell}=q_{\ell\ell}^{-1}=q_{jj}^{3}$,  $\qti_{jk}=q_{jj}^{-1}$ and $\qti_{ik}=\qti_{i\ell}=\qti_{j\ell}=1$;
\item \label{item:[[xijkl,xj],xk]-b}
$q_{ii}=\qti_{ij}^{\ -1} = -q_{\ell \ell}^{-1} = -\qti_{kl}$, $q_{jj}=\qti_{jk}=q_{kk}=-1$ and $\qti_{ik}=\qti_{i\ell}=\qti_{j\ell}=1$;
\end{enumerate}
Then $[[x_{ijk\ell},x_j]_c,x_k]_c=q_{jk}(\qti_{ij}^{\,-1}-q_{jj})[[x_{ijk\ell},x_k]_c,x_j]_c$ in $\toba$.
\end{lemma}

\pf
Assume that $ x_\beta  :=  [[x_{ijk\ell},x_j]_c,x_k]_c-q_{jk}(\qti_{ij}^{\,-1}-q_{jj})[[x_{ijk\ell},x_k]_c,x_j]_c\ne 0$.

A set of vertices as in \ref{item:[[xijkl,xj],xk]-a} only appears when $\bq$ is of type $\superf$. The proof is analogous to that of \cite[Lemma 3.24]{ACSa}.

A set of vertices as in \ref{item:[[xijkl,xj],xk]-b} only appears when $\bq$ is of type $\bgl(4)$.
The proof is analogous to that of \cite[Lemma 3.10 (ii)]{ACSa2}.
\epf

\section{Eminent pre-Nichols algebras}\label{sec:Eminent-infinite}

In this section we determine the eminent pre-Nichols algebras for all matrices $\bq$ considered here. In most cases the Nichols algebra is eminent, see Theorem \ref{thm:nichols-is-eminent}. For each one of the exceptions we describe explicitly the pre-Nichols algebra by generators and relations.

\subsection{Nichols algebra with trivial poset $\prefd$}

Let $\bq$ be such that $\GK \toba_{\bq}<\infty$ and $\dim \toba_{\bq}=\infty$. The following result characterizes in which cases the poset $\prefd$ has a unique element, the Nichols algebra:

\begin{theorem}\label{thm:nichols-is-eminent}
Let $\bq$ be a braiding matrix such that $\GK \toba_{\bq}<\infty$, $\dim\toba_{\bq}=\infty$ and the Dynkin diagram of $\bq$ is connected. 
If $\bq$ is not of one of the following types:
\begin{itemize}[leftmargin=*]
\item $\superqa{3}{q}{\{2\}}$ or $\superqa{3}{q}{\{1,2,3\}}$ 
with $q\not\in \G_{\infty}$,
\item $D(2,1;\alpha)$ where at least 2 of the 3 parameters are not roots of unity,
\end{itemize}
then $\toba_{\bq}$ is a eminent pre-Nichols algebra of $\bq$. Therefore $\toba_{\bq}$ is the only finite $GK$-dimensional pre-Nichols algebra of $\bq$.
\end{theorem}
\pf
The proof is analogous to that of \cite[Theorem 3.1]{ACSa2}, also using \cite[Remark 3.2]{ACSa2} together with the corresponding lemmas in $\S$ \ref{sec:lemmas}.
\epf

\subsection{Eminent pre-Nichols algebras of type Super A, exceptional cases.}

If $\bq$ is of type $\superqa{3}{q}{\{2\}}$ or $\superqa{3}{q}{\{1,2,3\}}$ with $q\not\in \G_{\infty}$, then we obtain analogous results to that of the corresponding cases when $\dim\toba_{\bq}<\infty$:

\begin{prop}\label{prop:superA-first-case-eminent}
If $\bq$ is of type $\superqa{3}{q}{\{2\}}$, then the quotient
\begin{align}\label{eq:superA-first-case-eminent}
\htoba_{\bq}=T(V)/\langle x_2^2,x_{13},x_{112},x_{332} \rangle,
\end{align}
is a eminent pre-Nichols algebra of $\bq$, with basis
\begin{align*}
B=\big\{x_3^ax_{23}^bx_2^cx_{12^23}^dx_{123}^ex_{12}^fx_1^g: \, b,c,e,f \in \{0,1\}, \, a,d,g\in \N_{0} \big\}
\end{align*}
where $x_{12^23}=[x_{123},x_2]_c$. Thus $\GK \htoba_{\bq} = 3$.
\end{prop}
\pf
Analogous to \cite[Proposition 5.5]{ACSa}.
\epf

\begin{remark}\label{rem:caseA1}
Here, $z=x_{12^23}$ is $q$-central and primitive in $\htoba_{\bq}$. Let $\Zc$ be the subalgebra of $\htoba_{\bq}$ generated by $z$. Then
$\Zc \hookrightarrow \htoba_{\bq} \twoheadrightarrow \toba_{\bq}$ is an extension of graded braided Hopf algebras.
\end{remark}

\begin{prop}\label{prop:superA-second-case-eminent}
If $\bq$ is of type $\superqa{3}{q}{\{1,2,3\}}$, then the quotient
\begin{align}\label{eq:superA-second-case-eminent}
\htoba_{\bq}=T(V)/\langle x_1^2, x_2^2, x_3^2, x_{213}, [x_{123},x_2]_c \rangle.
\end{align}
is a eminent pre-Nichols algebra of $\bq$, with basis
\begin{align*}
B=\big\{x_3^ax_{23}^bx_2^cx_{13}^dx_{123}^ex_{12}^fx_1^g: \, a,c,e,g\in \{0,1\}, \, b,d,f \in \N_{0} \big\}.
\end{align*}
Thus $\GK \htoba_{\bq} = 3$.
\end{prop}
\pf
Analogous to \cite[Proposition 5.9]{ACSa}.
\epf

\begin{remark}\label{rem:caseA2}
Here, $z=x_{213}$ is $q$-central and primitive in $\htoba_{\bq}$. Let $\Zc$ be the subalgebra of $\htoba_{\bq}$ generated by $z$. Then
$\Zc \hookrightarrow \htoba_{\bq} \twoheadrightarrow \toba_{\bq}$ is an extension of graded braided Hopf algebras.
\end{remark}

\subsection{Eminent pre-Nichols algebras of type $D(2,1;\alpha)$, exceptional cases.}\label{subsec:Eminent-D21a} The diagrams of type $D(2,1;\alpha)$ depends on three parameters $q,r,s$ such that $qrs=1$.
Let $\bq$ be a matrix of type $D(2,1;\alpha)$ such that exactly two of their parameters $q,r,s$ are not roots of one. 
Namely, set $M=\ord q$, $N=\ord r$, $L=\ord s$. We deal with the following cases:
\begin{align}\label{diagram:first-exceptional-case-D21a}
\xymatrix@C=30pt{\overset{q}{\circ} \ar@{-}[r]^{q^{-1}} &\overset{-1}{\circ}\ar@{-}[r]^{r^{-1}}\ &\overset{r}{\circ}}&&M&<\infty\quad N,L=\infty
\\
\label{diagram:second-exceptional-case-D21a}
\xymatrix@C=30pt{\overset{q}{\circ} \ar@{-}[r]^{q^{-1}} &\overset{-1}{\circ}\ar@{-}[r]^{r^{-1}}\ &\overset{r}{\circ}}&&L&<\infty\quad M,N=\infty
\\
\label{diagram:third-exceptional-case-D21a}
\xymatrix@C=30pt{
&\overset{-1}{\underset{1}{\circ}}\ar@{-}[ld]_{q}\ar@{-}[rd]^{r} & 
\\
\overset{-1}{\underset{2}{\circ}} \ar@{-}[rr]^{s} & &\overset{-1}{\underset{3}{\circ}}}&&M&<\infty\quad N,L=\infty
\end{align}


In the remainder part of this subsection we will introduce a pre-Nichols algebra and prove that it is the corresponding eminent one for each diagram.

\subsubsection{Type $D(2,1;\alpha)$ with Dynkin diagram (\ref{diagram:first-exceptional-case-D21a}) }\label{ssec:D21a-first-exception}
In this case, the Nichols algebra $\toba_{\bq}$ has the following presentation
\begin{align*}
x_2^2&=0, & x_{13}&=0, & x_{112}&=0, & x_{332}&=0, & x_1^M&=0,
\end{align*}

\begin{prop}\label{prop:D(2,1,a)-first-case}
Let $\bq$ be of type $D(2,1;\alpha)$ with Dynkin diagram (\ref{diagram:first-exceptional-case-D21a}). Then the algebra 
\begin{align}\label{eq:D(2,1,a)-first-case}
\htoba_{\bq}=T(V)/\langle x_2^2,x_{13},x_{112},x_{332}\rangle    
\end{align}
is an eminent pre-Nichols of $\bq$ with basis
\begin{align*}
B=\{x_3^{n_1}x_{23}^{n_2}x_2^{n_3}x_{12^23}^{n_{4}}x_{123}^{n_5}x_{12}^{n_6}x_1^{n_7}:n_2,n_3,n_5,n_6\in\{0,1\}\},
\end{align*}
where $x_{12^23}=[x_{123},x_2]_c$, and $\GK \htoba_{\bq}=3$. Also, if $\Zc$ is the algebra spanned by $x_1^M$ then $\Zc \hookrightarrow \htoba_{\bq} \twoheadrightarrow \toba_{\bq}$ is an extension of graded braided Hopf algebras.
\end{prop}
\pf
Let $I$ be the defining ideal of $\htoba_{\bq}$. By \cite[Remark 3.2]{ACSa2} (taking $J$ as the Hopf ideal that defines $\toba_{\bq}$), each generator of $I$ is primitive, so $\htoba_{\bq}$ is a pre-Nichols algebra. 

Let $\toba$ be a pre-Nichols algebra of $\bq$ such that $\GK\toba<\infty$. By Lemmas \ref{lem:mij=0}, \ref{lem:mij>0},\ref{lem:xi-no-Cartan} every defining relation of $\htoba_{\bq}$ vanishes in $\toba$. Then the canonical projection $T(V)\twoheadrightarrow \toba$ induces a surjective Hopf algebra map $\htoba_{\bq}\twoheadrightarrow \toba$. 
Thus, to show that $\htoba_{\bq}$, it only remains to prove that $\GK \htoba_{\bq}<\infty$.

By the proof of the \cite[Theorem 3.1]{An-crelle}, $B$ is a basis of $\htoba_{\bq}$. Let $\Zc'=\htoba^{\co \pi}$, where $\pi:\htoba_{\bq} \twoheadrightarrow \toba_{\bq}$ is the canonical projection. Since $\Zc\subseteq \Zc'$, from \cite[Lemma 2.4]{ACSa} we get
\begin{align*}
\mathcal{H}_{\htoba_{\bq}}&=\mathcal{H}_{\Zc'}\mathcal{H}_{\toba_{\bq}}\geq \mathcal{H}_{\Zc}\mathcal{H}_{\toba_{\bq}}\geq\\
&\geq \frac{1}{1-t_1^M}\frac{(1-t_1^M)(1+t_1t_2)(1+t_1t_2t_3)(1+t_2)(1+t_2t_3)}{(1-t_1)(1-t_1t_2^2t_3)(1-t_3)}.
\end{align*}
On the other hand $B$ is a basis of $\htoba_{\bq}$, so
\begin{align*}
\mathcal{H}_{\htoba_{\bq}}= \frac{(1+t_1t_2)(1+t_1t_2t_3)(1+t_2)(1+t_2t_3)}{(1-t_1)(1-t_1t_2^2t_3)(1-t_3)}
\end{align*}
so $\Zc=\Zc'$, $\Zc \hookrightarrow \htoba_{\bq} \twoheadrightarrow \toba_{\bq}$ is an extension of graded braided Hopf algebras and $\GK \htoba_{\bq}=3$.
\epf

\begin{remark}\label{rem:caseD1}
Here, $z=x_{1}^M$ is $q$-central and primitive in $\htoba_{\bq}$.
\end{remark}

\subsubsection{Type $D(2,1;\alpha)$ with Dynkin diagram \ref{diagram:second-exceptional-case-D21a} }\label{sssec:D21a-second-exception}

Let $\bq$ be of type $D(2,1;\alpha)$ with Dynkin diagram (\ref{diagram:second-exceptional-case-D21a}). The Nichols algebra $\toba_{\bq}$ has the following presentation

\begin{align*}
x_{13}&=0, & x_{2}^2&=0, & x_{112}&=0, & x_{332}&=0, & x_{12^23}^L&=0,
\end{align*}
where $x_{12^23}=[x_{123},x_2]_c$.

\begin{prop}\label{prop:D(2,1,a)-second-case}
Let $\bq$ be of type $D(2,1;\alpha)$ with Dynkin diagram (\ref{diagram:second-exceptional-case-D21a}). Then the algebra 
\begin{align}\label{eq:D(2,1,a)-second-case}
\htoba_{\bq}=T(V)/\langle x_2^2,x_{13},x_{112},x_{332}\rangle    
\end{align}
is an eminent pre-Nichols of $\bq$ with basis
\begin{align*}
B=\{x_3^{n_1}x_{23}^{n_2}x_2^{n_3}x_{12^23}^{n_{4}}x_{123}^{n_5}x_{12}^{n_6}x_1^{n_7}:n_2,n_3,n_5,n_6\in\{0,1\}\}
\end{align*}
so $\GK \htoba_{\bq} =3$. Aslo, if $\Zc$ be the algebra spanned by $x_{12^23}^M$ then $\Zc \hookrightarrow \htoba_{\bq} \twoheadrightarrow \toba_{\bq}$ is a degree-preserving extension of braided Hopf algebras.

\end{prop}
\pf
Let $I$ be the defining ideal of $\htoba_{\bq}$. By \cite[Remark 3.2]{ACSa2} (taking $J$ as the Hopf ideal that defines $\toba_{\bq}$), each generator of $I$ is primitive, so $\htoba_{\bq}$ is a pre-Nichols algebra. 

Let $\toba$ be a pre-Nichols algebra of $\bq$ such that $\GK\toba<\infty$. By Lemmas \ref{lem:mij=0}, \ref{lem:mij>0},\ref{lem:xi-no-Cartan} every defining relation of $\htoba_{\bq}$ vanishes in $\toba$. Then the canonical projection $T(V)\twoheadrightarrow \toba$ induces a surjective Hopf algebra map $\htoba_{\bq}\twoheadrightarrow \toba$. 
Thus, to show that $\htoba_{\bq}$, it only remains to prove that $\GK \htoba_{\bq}<\infty$.

By the proof of the \cite[Theorem 3.1]{An-crelle}, $B$ is a basis of $\htoba_{\bq}$. Let $\Zc'=\htoba^{\co \pi}$, where $\pi:\htoba_{\bq} \twoheadrightarrow \toba_{\bq}$ is the canonical projection. Since $\Zc\subseteq \Zc'$, from \cite[Lemma 2.4]{ACSa} we get

\begin{align*}
\mathcal{H}_{\htoba_{\bq}}&=\mathcal{H}_{\Zc'}\mathcal{H}_{\toba_{\bq}}\geq \mathcal{H}_{\Zc}\mathcal{H}_{\toba_{\bq}}\geq\\
&\geq \frac{1}{1-(t_1t_2^2t_3)^M}\frac{(1-(t_1t_2^2t_3)^M)(1+t_1t_2)(1+t_1t_2t_3)(1+t_2)(1+t_2t_3)}{(1-t_1)(1-t_1t_2^2t_3)(1-t_3)}.
\end{align*}
On the other hand $B$ is a basis of $\htoba_{\bq}$, so
\begin{align*}
\mathcal{H}_{\htoba_{\bq}}= \frac{(1+t_1t_2)(1+t_1t_2t_3)(1+t_2)(1+t_2t_3)}{(1-t_1)(1-t_1t_2^2t_3)(1-t_3)}
\end{align*}
so $\Zc=\Zc'$, $\Zc \hookrightarrow \htoba_{\bq} \twoheadrightarrow \toba_{\bq}$ is an extension of graded braided Hopf algebras and $\GK \htoba_{\bq}=3$.


\epf
\begin{remark}\label{rem:caseD2}
$z=x_{12^23}^L$ is $q$-central and primitive in $\htoba_{\bq}$.
\end{remark}

\subsubsection{Type $D(2,1;\alpha)$ with Dynkin diagram \ref{diagram:third-exceptional-case-D21a} }\label{sssec:D21a-third-exception}

Let $\bq$ be of type $D(2,1;\alpha)$ with Dynkin diagram (\ref{diagram:third-exceptional-case-D21a}). The Nichols algebra $\toba_{\bq}$ has the following presentation

\begin{align*}
x_2^2&=0, & x_{2}^2&=0, & x_{3}^2&=0, & x_{12}^M&=0,
\end{align*}
\begin{align*}
x_{123}= q_{12}(1-s)x_{2}x_{13}-\frac{1-s}{q_{32}(1-r)}[x_{13},x_{2}]_c.
\end{align*}

\begin{prop}\label{prop:D(2,1,a)-third-case}
Let $\bq$ be of type $D(2,1;\alpha)$ with Dynkin diagram (\ref{diagram:third-exceptional-case-D21a}). The algebra $\htoba_{\bq}$ defined by the following relations
\begin{align}\label{eq:D(2,1,a)-third-case}
\begin{aligned}
x_1^2&=0, \qquad x_{2}^2=0, \qquad x_{3}^2=0,
\\
x_{123} &= q_{12}(1-s)x_{2}x_{13}-\frac{1-s}{q_{32}(1-r)}[x_{13},x_{2}]_c.
\end{aligned}
\end{align}

is an eminent pre-Nichols of $\bq$, with basis
\begin{align*}
B=\{x_3^{n_1}x_{23}^{n_2}x_2^{n_3}x_{123}^{n_4}x_{13}^{n_5}x_{12}^{n_6}x_1^{n_7}:n_1,n_3,n_4,n_7\in\{0,1\}\}
\end{align*}
so $\GK \htoba_{\bq} =3$. Aslo, if $\Zc$ be the algebra spanned by $x_{12}^M$ then $\Zc \hookrightarrow \htoba_{\bq} \twoheadrightarrow \toba_{\bq}$ is a degree-preserving extension of braided Hopf algebras.
\end{prop}
\pf
Let $I$ be the defining ideal of $\htoba_{\bq}$. By \cite[Remark 3.2]{ACSa2} (taking $J$ as the Hopf ideal that defines $\toba_{\bq}$), each generator of $I$ is primitive, so $\htoba_{\bq}$ is a pre-Nichols algebra. 

Let $\toba$ be a pre-Nichols algebra of $\bq$ such that $\GK\toba<\infty$. By Lemmas \ref{lem:xi-no-Cartan}, \ref{lem:triangle-x_{ijk}} every defining relation of $\htoba_{\bq}$ vanishes in $\toba$. Then the canonical projection $T(V)\twoheadrightarrow \toba$ induces a surjective Hopf algebra map $\htoba_{\bq}\twoheadrightarrow \toba$. 
Thus, to show that $\htoba_{\bq}$, it only remains to prove that $\GK \htoba_{\bq}<\infty$.

 By the proof of the \cite[Theorem 3.1]{An-crelle}, $B$ is a basis of $\htoba_{\bq}$. Let $\Zc'=\htoba^{\co \pi}$, where $\pi:\htoba_{\bq} \twoheadrightarrow \toba_{\bq}$ is the canonical projection. Since $\Zc\subseteq \Zc'$, from \cite[Lemma 2.4]{ACSa} we get
\begin{align*}
\mathcal{H}_{\htoba_{\bq}}&=\mathcal{H}_{\Zc'}\mathcal{H}_{\toba_{\bq}}\geq \mathcal{H}_{\Zc}\mathcal{H}_{\toba_{\bq}}\geq\\
&\geq \frac{1}{1-(t_1t_2)^M}\frac{(1-(t_1t_2)^M)(1+t_1)(1+t_1t_2t_3)(1+t_2)(1+t_3)}{(1-t_1t_2)(1-t_1t_3)(1-t_2t_3)}.
\end{align*}
On the other hand $B$ is a basis of $\htoba_{\bq}$, so
\begin{align*}
\mathcal{H}_{\htoba_{\bq}}= \frac{(1+t_1)(1+t_1t_2t_3)(1+t_2)(1+t_3)}{(1-t_1t_2)(1-t_1t_3)(1-t_2t_3)}
\end{align*}
so $\Zc=\Zc'$, $\Zc \hookrightarrow \htoba_{\bq} \twoheadrightarrow \toba_{\bq}$ is an extension of graded braided Hopf algebras and $\GK \htoba_{\bq}=3$.
\epf
\begin{remark}\label{rem:caseD3}
$z=x_{12}^M$  is $q$-central and primitive in $\htoba_{\bq}$.
\end{remark}

\subsection{Nichols algebras with not connected Dynkin diagram}
Let $\bq$ be such that $\GK \toba_{\bq}<\infty$, $\dim \toba_{\bq}=\infty$ and the Dynkin diagram of $\toba_{\bq}$ satisfy
We consider $\I=\cup_{d\in C}\I^{(d)}$ such that the $\I^{(d)}=\{j_{(1)},j_{(2)},...,j_{(k_d)}\}$ as a set of vertices of connected component of diagram of $\bq$. Thus $V_d=\langle x_j|j\in \I^{(d)}\rangle$ are a braided subspace of $V$ and we denote $\bq_d$ a the braided matrix of $V_d$.\\

The following remark follows from the rank 3 classification and will be useful to prove the main theorem of this section.
\begin{remark}
    If $V=\ku_{1}+\ku_{2}+\ku_{3}$ a vector braided space with braided matrix $\bq$ such that $\GK \toba_{\bq}<\infty$, $\dim \toba_{\bq}=\infty$ and the Dynkin diagram of $\bq$ is a 3-cycle, then it satisfies the following conditions:
    \begin{itemize}
        \item At least two of the vertices labels are -1.
        \item If the label of a vertex is not -1, then that label is not a root of unity.
        \item If $q_{11}\ne -1$ then $\qti_{12}=\qti_{13}=q_{11}^{-1}$ or $\{\qti_{13},\qti_{12}\}=\{q_{11}^{-1},q_{11}^{-2}\}$.
    \end{itemize}
\end{remark}

\begin{theorem}\label{thm:eminent-not-connected-diagram}
Let $\bq$ be a braided matrix such that $\GK \toba_{\bq}<\infty$ and satisfy the following conditions
\begin{itemize}
    \item For all $d\in C$, either $|k_d|\geq 2$ or else  $|k_d|=1$ with $q_{ii}\ne 1$, $i\in \I^{(d)}$,
    \item $\dim \toba_{\bq_d}=\infty$ for all $d\in C$,
\end{itemize}

Then the pre-Nichols algebra $\htoba_{\bq}$ defined by generators $x_i,i\in \I$ and the following relations:
\begin{enumerate}[leftmargin=*,label=\rm{(\alph*)}]
    \item\label{itm-a-thm-non-connected-relations} defining relations of $\htoba_{\bq_d}$, $d\in C$,
    \item\label{itm-b-thm-non-connected-relations} $x_{ij}=0$ if $i\in\I^{(a)}$, $j\in\I^{(b)}$, $a\neq b$.
\end{enumerate}
is a eminent pre-Nichols algebra of $\bq$ and is isomorphic to $\underline{\bigotimes}_{d\in C}{\htoba_{\bq_{d}}}$.
\end{theorem}
\pf
Let $\toba$ a finite GK-dimensional pre-Nichols algebra of $\toba$. Note that by Lemmas in $\S 3$ \cite[\S 3]{ACSa} and \cite[\S 3.2]{ACSa2},  the relations in \ref{itm-a-thm-non-connected-relations} are satisfied in $\toba$, in efect, their supports are contained in $\I^{(a)}$ for some $a$ and the Dynkin diagram of the braided vector subspace $V_a$ is connected.

To prove that relations \ref{itm-b-thm-non-connected-relations} are satisfied in $\toba$ we proceed analogously to $\S 3$. Let $i,j\in \I$ be as in \ref{itm-b-thm-non-connected-relations}, we have to $x_{ij}\in\Pc(\toba)$. Suppose that $x_{ij}\neq 0$ in $\toba$, then the Dynkin diagram of braided vector space $W:=\ku x_{j}+\ku x_{j}+\ku x_{ij}$ is
\begin{align*}
\xymatrix@C=30pt{\overset{q_{ii}}{\underset{i}{\circ}}\ar@{-}[rr]^{q_{ii}^2} && \overset{q_{ii}^2q_{jj}^2}{\underset{ij}{\circ}}\ar@{-}[rr]^{q_{jj}^2} && \overset{q_{jj}}{\underset{j}{\circ}}
}
\end{align*}

Suppose that $|k_a|=|k_b|=1$, then $q_{ii},q_{jj}\not\in\G_{\infty}$ since $\dim \toba_{\bq_a}=\dim \toba_{\bq_b}=\infty$. Since $q_{ii}^2$ is a label of an edge, by rank 3 classification we have that $q_{ii}^2=-1$ or exist $m\in \{1,2,3\}$ such that $q_{ii}^mq_{ii}^2=1$. We get a contradiction since $q_{ii}\not\in \G_{\infty}$ and we conclude that $|k_a|>1$ or $|k_b|>1$. 

Suppose that $|k_a|>1$. Let $k\in \I^{(a)}$ be such that $\qti_{ik}\ne 1$, the Dynkin diagram of $W=\ku x_k+\ku x_{i}+\ku x_{j}+\ku x_{ij}$ is 
\begin{align*}
\xymatrix@C=30pt{&&&&\overset{q_{ii}}{\underset{i}{\circ}}\ar@{-}[dd]_{\qti_{ik}}\ar@{-}[dll]_{q_{ii}^2}\\ \overset{q_{jj}}{\underset{j}{\circ}}\ar@{-}[rr]^{q_{jj}^2}
&&\overset{q_{ii}^2q_{jj}^2}{\underset{ij}{\circ}}\ar@{-}[rrd]^{\qti_{ik}} &&\\
&&&&\overset{q_{kk}}{\underset{k}{\circ}}\\
}
\end{align*}
If $q_{ii}\ne -1$ then the Dynkin diagram of $\ku x_k+\ku x_{i}+\ku x_{ij}$ is a 3-cycle and by remark, $q_{ii}^2q_{jj}^2=q_{kk}=-1$. By inspection, the only possibility for $W$ is of type $\superf$, but it is impossible because in this case we have $q_{ii}=q\not\in\G_{\infty}$ and $\qti_{ik}=q^{-1}=q^2$. We get a contradiction, therefore $q_{ii}=-1$. The Dynkin diagram of $\ku x_j+\ku x_k+\ku x_{i}+\ku x_{ij}$ is

\begin{align*}
\xymatrix@C=30pt{\overset{q_{jj}}{\underset{j}{\circ}}\ar@{-}[rr]^{q_{jj}^2} && \overset{q_{jj}^2}{\underset{ij}{\circ}}\ar@{-}[rr]^{\qti_{ik}} && \overset{q_{kk}}{\underset{k}{\circ}}\ar@{-}[rr]^{\qti_{ik}} && \overset{-1}{\underset{i}{\circ}}
}
\end{align*}

Note that by rank 4 classification, each label in the above diagram satisfies that it is not a root of unity or is -1. It is concluded that $\qti_{ik},q_{jj},q_{jj}^2,q_{kk}\not\in \G_{\infty}-\G_{2}$. So, we have that $q_{jj}^2 \ne -1$ otherwise $q_{jj}\in\G_{4}-\G_2$. Once again, due to the rank 4 classification, we have that if an adjacent vertex and edge have equal labels, then it is of type $\bgl(4)$ and the mentioned labels are equal to $-1$. The latter is a contradiction since the vertex $ij$ is adjacent to the edge $\{j,ij\}$ and $q_{jj}^2\ne -1$.

Due to the previous contradiction we have that $x_{ij}=0$ and this concludes what we wanted to prove.
\epf

\begin{remark}
    If $\bq$ is as in Theorem \ref{thm:eminent-not-connected-diagram}, then we can compute the poset of graded pre-Nichols algebras with finite $\GK$. The elements of the poset would be parameterized by $\underline{\bigotimes}_{d\in C}{\toba_{d}}$ where $\toba_d=\toba_{\bq_d}$ or $\toba_d=\htoba_{\bq_d}$.
\end{remark}
\textbf{Acknowledgments}
I thank Iv\'an Angiono for suggesting this problem and for reading our work.





%

\end{document}